# Towards a classical proof of the 4-colour theorem

Peter Dörre

Fachhochschule Südwestfalen (University of Applied Sciences)
Frauenstuhlweg 31, D-58644 Iserlohn, Germany
Email: doerre.peter(at)fh-swf.de

## Abstract

It is easy to colour the subset of triangulations with minimum degree less than 5 by induction on the graph order *n*, whereas the case with minimum degree equal to 5 cannot be coloured by recolouring the neighbours of a single 5-valent vertex. Instead of trying to colour appropriate subgraphs composed of at least 5-valent vertices, the new strategy is to use the colouring information already provided by induction in a constructive way.

## 1 Introduction

The 4-colour theorem was first proved by Appel and Haken (see e.g. [1]) in 1976 with the help of extensive computer calculations. Improved and independent versions of this type of proof exist since 1996 [2]. We note that in all these works, Kempe chains are used as a proof-tool, originally introduced by Kempe [3] in his attempted proof. Note that Kempe's proof strategy is correct, as long as in the induction step, a vertex of degree less than 5 can be found: then there is no need to recolour a 5-valent vertex.

Apart from some 4-cycles and 5-cycles, all graphs considered here are triangulations. Before we begin, we subdivide the set of graphs into two mutually exclusive subsets. Let $G_<$ denote a graph from the subset with $\delta(G) < 5$, and $G_=$ a graph from the subset with $\delta(G) = 5$, respectively.

Assume we have coloured all graphs with $\delta(G) < 5$ by induction on the graph order *n* along the lines of Kempe's proof, the question arises whether it is possible to construct from this subset the omitted graphs with $\delta(G) = 5$ with (proper) 4-colourings. Our strategy is to use the colourings induction found for us, and to complete them.

## 2 Theorem

*The chromatic number of a planar graph is at most four.*

## 3 Proof of the theorem

### 3.1 Graphs with minimum degree less than 5

In a proof by induction on the graph order *n,* when we carry out the induction step, the idea is to find a suitable uncoloured vertex *v*, to remove it from the graph *G*, and to colour the remaining graph $G - v$ by the induction hypothesis. If we manage to obtain by recolouring a 3-colouring of the cycle of neighbours of the removed vertex *v*, then we may insert this vertex again, coloured with the unused colour, to complete the colouring of *G*.

In the subset of graphs with $\delta(G) < 5$, we always find a vertex *v* with degree at most 4. Suppose $C_4(v) = v_1 \cdots v_4$ is such a cycle of neighbours of a 4-valent vertex *v*. Now we introduce a small but crucial modification: we remove *v* and add the edge $v_2 v_4$ (see Fig. 1), and obtain the graph $G_< - v + v_2 v_4$ to be coloured by the induction hypothesis. By using a 4-cycle with an edge $v_2 v_4$ in its interior, only a 3- or 4-colouring is obtained. Possible colourings are shown in Fig. 1.

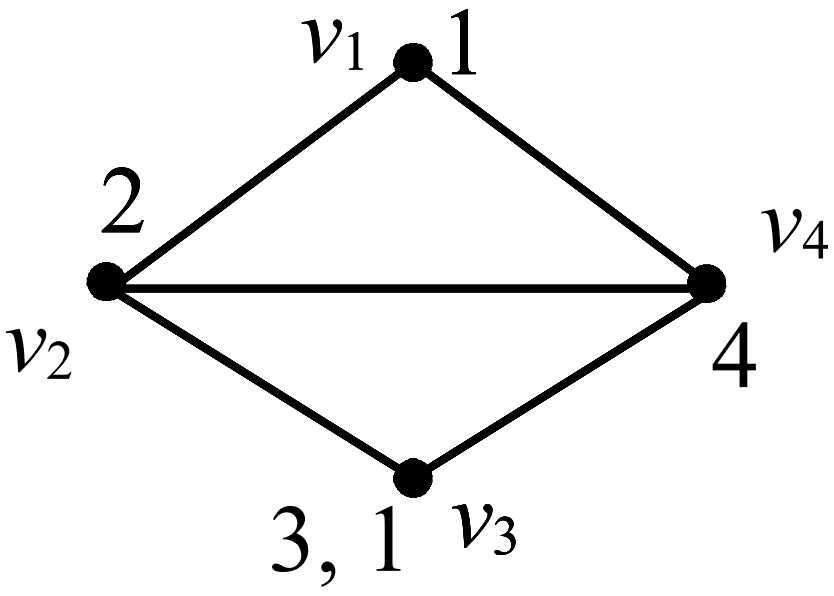

**Fig. 1**. The cycle $C_4 \subseteq G_<$ of neighbours of the 4-valent vertex *v*

If $C_4(v)$ is not yet 3-coloured, assume a colour sequence 1-2-3-4 in counterclockwise order. There cannot exist two Kempe chains $P_{13}(v_1, v_3)$ (between $v_1$ and $v_3$) and $P_{24}(v_2, v_4)$ (between $v_2$ and $v_4$) simultaneously: at least one must be absent. Hence we can always recolour in one endpoint of the non-existing Kempe chain to obtain a 3-colouring on $C_4(v)$, and colour *v* with the remaining colour.

Note that for recolouring, the edge $v_2 v_4$ is no longer relevant and can be removed, to be replaced by the vertex *v* connected to his neighbours, and coloured with the colour $c(v)$ not on $C_4$, to complete the colouring of $G_<$.

To prepare for another evaluation of the induction step, consider the following 5-cycle $C_5 \subseteq G_<$ and its interior (see Fig. 2).

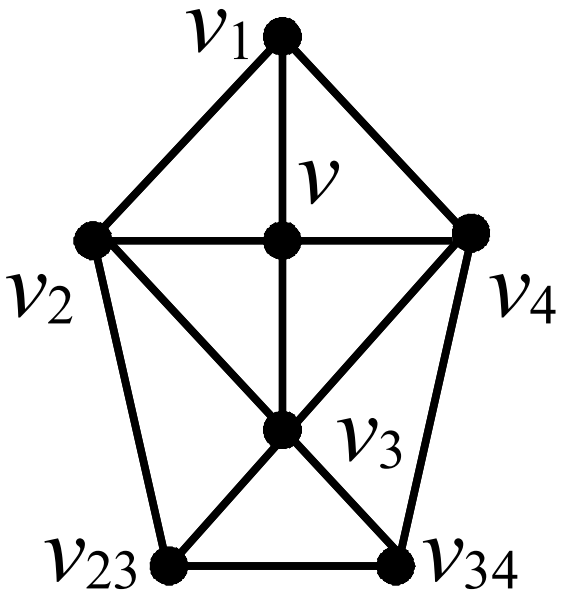

**Fig. 2**. Denotation of vertices of the cycle $C_5$ and its interior

Note that we denoted the vertices on $C_5$ that are not on $C_4(v)$ by their adjacency: $v_{23}$ is adjacent to the pair $(v_2, v_3)$, and $v_{34}$ to $(v_3, v_4)$.

When we created $G_< - v + v_2 v_4$, vertex $v_3$ changes into a 4-valent vertex, such that every graph which contains the graph $C_5{}^*$ (see Fig. 3) belongs to the subset with $\delta(G) < 5$ and hence can be coloured by the induction hypothesis. Note that we are no longer interested in reinserting vertex *v*.

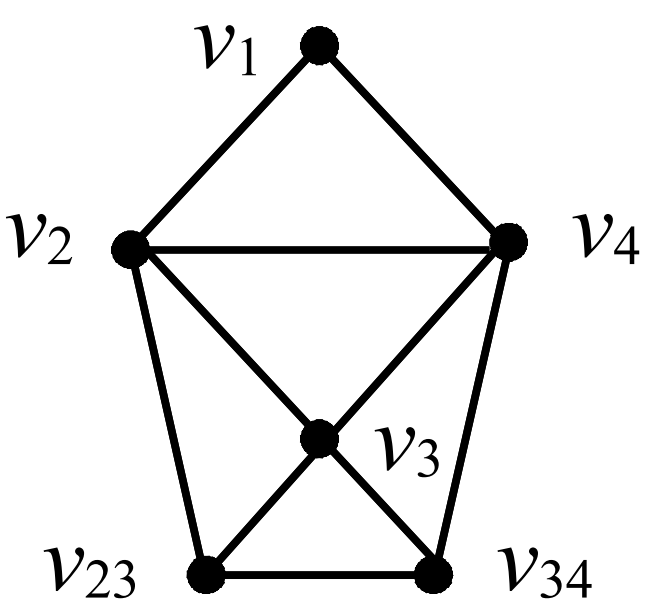

**Fig. 3**. The graph $C_5{}^* \subseteq G_<$

**Selection**. With respect to the subset of graphs with $\delta(G) < 5$, we now focus our exclusive interest on those special graphs, which contain $C_5{}^*$ as a subgraph, but only once, and no other vertices of degree 3 or 4. These special graphs are already coloured.

Next we investigate the colourings obtained for $C_5{}^* \subseteq G_<$ by the induction hypothesis in detail. Note that these colourings are valid for every graph $G_< - v + v_2 v_4$.

Our goal is to find a 3-colouring of the 5-cycle $C_5 = v_1v_2v_{23}v_{34}v_4$, such that after removing its interior (with the last 4-valent vertex), we can insert a 5-valent vertex $v'$ coloured with the unused colour.

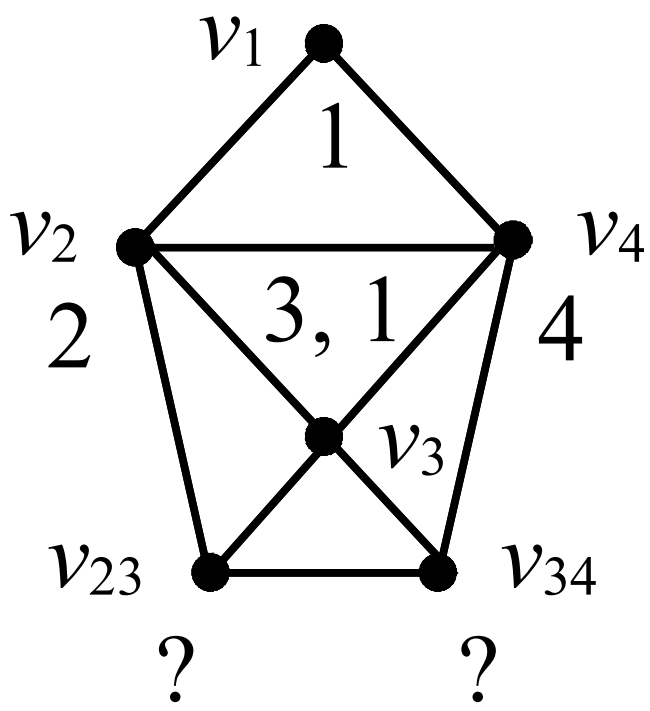

**Fig. 4**. Colourings obtained for $C_5{}^* \subseteq G_<$

### 3.2 Construction of colourings for graphs with minimum degree equal to 5

**Part 1**. The cycle $C_4$ received a 4-colouring by the induction hypothesis (Fig. 5)

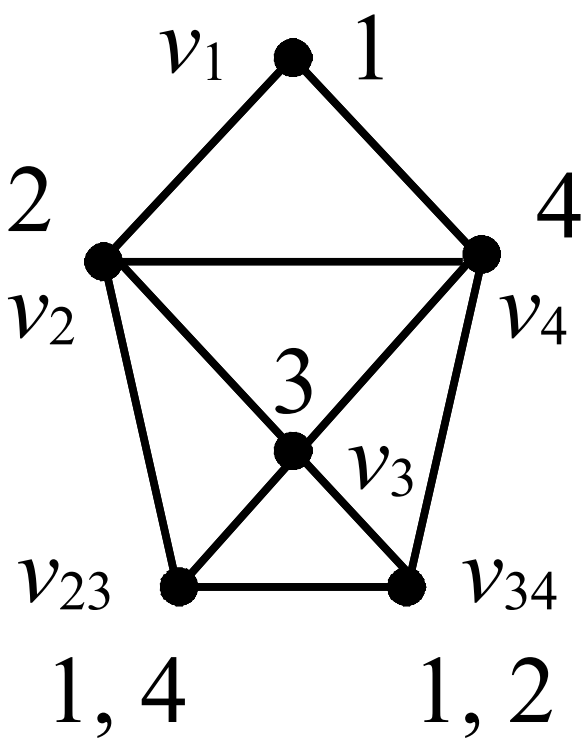

**Fig. 5**. The graph $C_5{}^*$ with a 4-colouring of $C_4$

As shown in Fig. 5, the pair of vertices $(v_{23}, v_{34})$ comes out with the following pairs of colours: (1, 2), (4, 1), and (4, 2), which provide for $C_5$ the three colour sequences 1-2-1-2-4, 1-2-4-1-4, and 1-2-4-2-4. As all 3 sequences describe 3-colourings, we complete the colouring of $G_=$ with $c(v') = 3$ for the 5-valent vertex $v'$ in the interior of $C_5$. Note that here the existence or non-existence of Kempe chains plays no role in the construction process.

**Part 2**. The cycle $C_4$ received a 3-colouring by the induction hypothesis (Fig. 6)

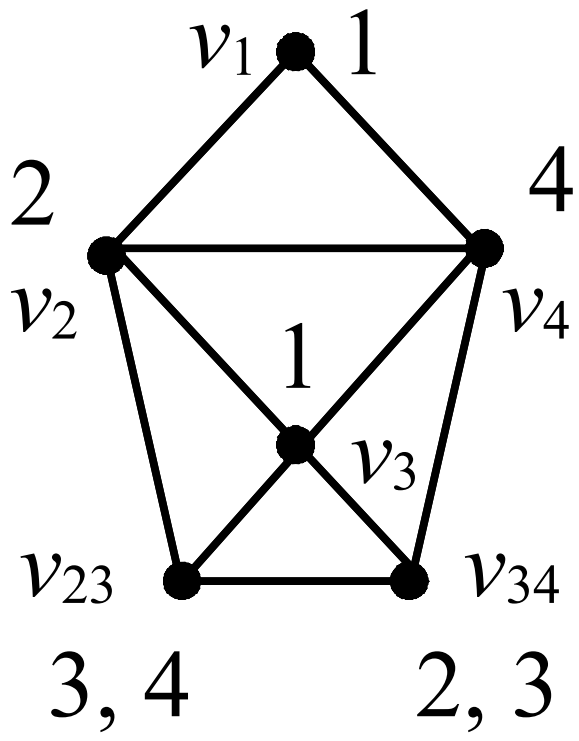

**Fig. 6**. The graph $C_5$ * with a 3-colouring of $C_4$

Now we obtain one 3-colouring, namely the colour sequence 1-2-4-2-4 (see Fig. 7), and two 4-colourings 1-2-3-2-4 and 1-2-4-3-4, which cannot be recoloured to a 3-colouring on $C_5$.

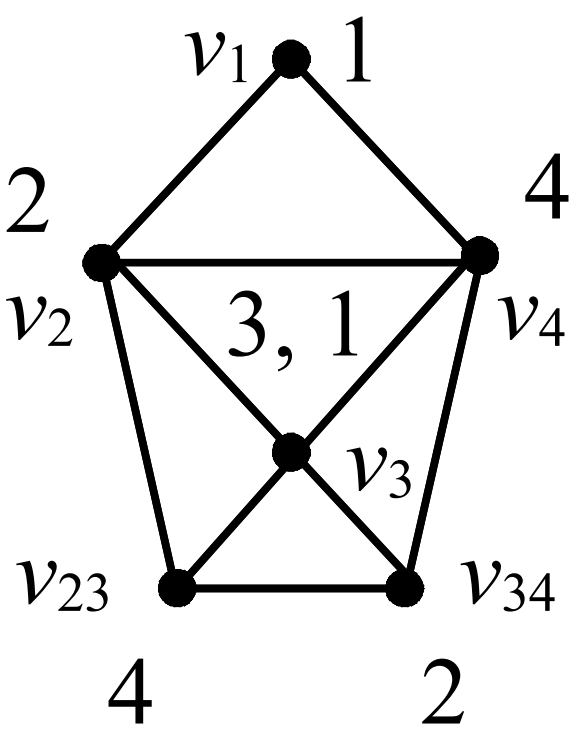

**Fig. 7**. Common colour sequence 1-2-4-2-4 (Part 1 and 2)

Note that the colour sequence 1-2-4-2-4 also occurs in the colourings proposed by the induction hypothesis in Part 1.

**Conclusion**. The theorem implies that for a graph $G_=$, there is at least one 4-colouring. As a result, the induction hypothesis provided us in all cases with one colour sequence of 3 colours, namely 1-2-4-2-4, which allows us to complete the colouring of $G_=$ by the choice $c(v') = 3$. Note that this colouring of the 5-cycle of neighbours of $v'$ shows a short Kempe chain $P_{24}(v_2, v_4)$ along $C_5$, thus excluding a chain $P_{13}(v_1, v_3)$, irrespective of the colour of vertex $v_3$ (3 or 1). Furthermore, this result describes the general situation at the 5-cycles around coloured 5-valent vertices: one single colour, and two double colours alternatingly distributed on $C_5$, thus defining a short Kempe chain of the two double colours.

**Summary**. Planar graphs with minimum degree less than 5 are coloured by induction on the graph order *n*, the other planar graphs with minimum degree equal to 5 are coloured by construction, using the colouring information provided by the induction hypothesis.

## Acknowledgments

This paper is written in honour of A. B. Kempe who provided the proof-tool.
I am indebted to Douglas B. West and Reinhard Diestel for helpful comments on previous versions of the paper.